\documentclass[10pt,a4paper,twoside]{article}

\usepackage{amsfonts, amssymb, amsmath, amsthm}
\usepackage{mathrsfs} % \mathscr
\usepackage{latexsym}%utilisation des symboles LaTeX pour avoir un beau LaTeX
\usepackage{enumerate}
\usepackage{multicol}
\usepackage{verbatim}

\usepackage{url}
\usepackage{csquotes}
\usepackage[colorlinks=true]{hyperref}
\hypersetup{citecolor=blue, linkcolor=blue}

\newtheorem{thm}{Theorem}%[section]

\let\ve=\varepsilon

\title{A probabilistic free-proof of an explicit Croot-\L aba-Sisask Lemma%
}
\author{Olivier Ramar\'e}

\begin{document}
% \address[O. Ramar\'e]{CNRS/ Institut de Math\'ematiques de Marseille, Aix 
% Marseille Universit\'e, U.M.R. 7373, Site Sud, Campus de Luminy, Case 907, 
% 13288 
% Marseille Cedex 9, France.}
% \email{olivier.ramare@univ-amu.fr}

% %\date{\sl January, the 7th of 2004}
% \subjclass[2010]{Primary: , Secondary: }

% \keywords{Class Field}

\maketitle

\begin{abstract}
  %\texttt{File \jobname.tex} 
  This note proposes a probabilistic language-free proof of the
  famous Croot-\L aba-Sisask Lemma. In between, we do the same for the
  Khintchine and Marcinkiewicz-Zygmund inequalities and explicitate the
  implied constants.
\end{abstract}

%{\small \tableofcontents}

%%%%%%%%%%%%%%%%%%%%%%%%%%%%%%%%%%%%%%%%%
\section{Introduction}
%%%%%%%%%%%%%%%%%%%%%%%%%%%%%%%%%%%%%%%%%

We propose a proof free of probabilistic language
of (a variant of) \cite[Lemma 3.2]{Croot-Laba-Sisask*13} by
E.~Croot, I.~\L aba and O. Sisask, which we now state.
%%%%%
We follow their notation and in particular, when $z$ is a complex number,
$z^\circ$ is
defined  to be $z/|z|$ when $z\neq0$, and to be~0
when $z=0$.
%%%%%%%%%%%%%%%%%%%%%%%%%%%%%%%%%%%%
\begin{thm}
  \label{CLS}
  Let $(X,\mu)$ be a probability space and $p\ge2$.
  Given a function $f$ in the form
  \begin{equation*}
    f=\sum_{k\le K}\lambda_k g_k
  \end{equation*}
  where $(g_k)_k$ is a collection of measurable functions on~$X$ of
  $L^p(\mu)$-norm at most~1.
  Let $\ve>0$. There exists an $L$-tuple
  $(k_1,\cdots,k_\ell)\in\{1,\cdots,K\}^L$ of length $L\le 20p/\ve^2$ such that
  \begin{equation*}
    \int_X\biggl|\frac{f(x)}{\|\lambda\|_1}
    -\frac{1}{L}\sum_{\ell\le L}\lambda_{k_\ell}^\circ g_{k_\ell}(x)\biggr|^p d\mu\le \ve^p,
  \end{equation*}
  where $\|\lambda\|_1=\sum_{k\le K}|\lambda_k|$.
\end{thm}
%%%%%%%%%%%%%%%%%%%%%%%%%%%%%%%%%%%%
This theorem has its origin in the paper
\cite{Croot-Sisask*10} by E.~Croot and O.~Sisask. We refer to the
paper~\cite{Pierce*21} by L.~Pierce for much deeper background on the
Khintchine and Marcinkiewicz-Zygmund inequalities.

%%%%%%%%%%%%%%%%%%%%%%%%%%%%%%%%%%%%%%%%%
\section{An upper explicit Khintchine Inequality}
%%%%%%%%%%%%%%%%%%%%%%%%%%%%%%%%%%%%%%%%%

%%%%%%%%%%%%%
\begin{thm}
  \label{KI}
  We have, when $p\ge1$,
  \begin{equation*}
    (1/2^N)\sum_{(\ve_n)\in\{\pm1\}^N}\biggl|\sum_{n\le
      N}c_n\ve_n\biggr|^p
    \le
    p^{p/2}\biggl(\sum_{n\le N}|c_n|^2\biggr)^{p/2}.
  \end{equation*}
\end{thm}
%%%%%%%%%%%%%
This is only half of the Khintchine Inequality and in a special
context, but this will be enough for us. We followed
\cite[Chapter 10, Theorem 1, page 354]{Show-Teicher*97} by Y.S.~Show
and H.~Teicher.

%%%%%%%%%%%%%%%%
\begin{proof}
  Let us start with $p=2k\ge2$, so that we may open the inner sum and
  get
  \begin{align*}
    2^NS(2k)
    &=
    \sum_{(\ve_n)\in\{\pm1\}^N}\biggl|\sum_{n\le
      N}c_n\ve_n\biggr|^p
    \\&=
    \sum_{\substack{s_1+s_2+\cdots+s_{N}=2k,\\
        s_n\ge0}}
    \binom{2k}{s_1,s_2,\cdots,s_N}
    \prod_{1\le n\le N}c_n^{s_n}
    \sum_{(\ve_n)\in\{\pm1\}^N}
    \prod_{n\le N}\ve_n^{s_n}
  \end{align*}
  by the multinomial theorem. The inner summand vanishes as soon as
  some $s_n$ is odd, whence, by letting $2t_n=s_n$, we get
  \begin{align*}
    2^NS(2k)
    &=
    \sum_{\substack{t_1+t_2+\cdots+t_{N}=k,\\
        t_n\ge0}}
    \binom{2k}{2t_1,2t_2,\cdots,2t_N}
    \prod_{n\le N}(c_n^2)^{t_n}
    \\&\le
    C
    \sum_{\substack{t_1+t_2+\cdots+t_{N}=k,\\
        t_n\ge0}}
    \binom{k}{t_1,t_2,\cdots,t_N}
    \prod_{n\le N}(c_n^2)^{t_n}
    =
    C \biggl(\sum_{n\le N}c_n^2\biggr)^k
  \end{align*}
  where
  \begin{align*}
    C
    &=
    \max\binom{2k}{2t_1,2t_2,\cdots,2t_N}\binom{k}{t_1,t_2,\cdots,t_N}^{-1}
    \\&\le
    \max\frac{2k(2k-1)\cdots (k+1)}{\prod_j 2t_j(2t_j-1)\cdots
    (t_j+1)}
    \\&\le\max\frac{k^k}{2^{t_1+\cdots+t_k}}\le (k/2)^k.
  \end{align*}
  As $S(p)$ is increasing, we simply choose $k=\lceil p/2\rceil$. This
  gives us $2k\ge p+2$ and thus
  \begin{equation*}
    (k/2)^k\le (p+2)^{1+p/2}\le (30\, p)^{p/2}.
  \end{equation*}
  This concludes the main part of the proof, except for the
  constant~30. We will not continue the proof but simply refer to the
  paper~\cite{Haagerup*81} by U.~Haagerup who shows that best
  constant is (be careful: the abstract of this paper misses a closing
  parenthesis for the value of $B_p$, but the value of $B_p$ displayed
  in the middle of page 232 misses a squareroot-sign around the $\pi$,
  as an inspection of the proof at the end the paper rapidly reveals)
  \begin{equation*}
    \begin{cases}
      1&\quad\text{when $0<p\le 2$},\\
      \sqrt{2}\biggl(\frac{\Gamma((p+1)/2)}{\sqrt{\pi}}\biggr)^{1/p}
      &\quad\text{when $2< p$}.
    \end{cases}
  \end{equation*}
  We readily check that this implies that the constant~1 rather
  than~30 is admissible.
  % f(p)= sqrt(2)*(gamma((p+1)/2)/sqrt(Pi))^(1/p)/sqrt(p)
  %Concerning the lower bound ...
\end{proof}
%%%%%%%%%%%%%%%%

%%%%%%%%%%%%%%%%%%%%%%%%%%%%%%%%%%%%%%%%%
\section{An upper explicit Marcinkiewicz-Zygmund Inequality}
%%%%%%%%%%%%%%%%%%%%%%%%%%%%%%%%%%%%%%%%%
%\cite[Chapter 10, Theorem 2, page 356]{Show-Teicher*97}

%%%%%%%%%%%%
\begin{thm}
  \label{MZI}
  Let $(X,\mu)$ be a probability space.
  When $p\ge1$, let $(f_n)_{n\le N}$ be a system of functions such that
  $\int_X f_n(x)d\mu=0$. We have
  \begin{multline*}
    \int_{(x_n)\in X^N}\biggl|\sum_{1\le n\le
      N}f_n(x_n)\biggr|^pd(x_n)
    \\\le
    (4p)^{p/2}
    \int_{(x_n)\in X^N}\biggl(\sum_{1\le n\le
      N}|f_n(x_n)|^2\biggr)^{p/2}d(x_n).
  \end{multline*}
\end{thm}
%%%%%%%%%%%%
The power of this inequality is that the implied constant do not
depend on~$N$, the effect of some orthogonality.
Again, this is only half of the Marcinkiewicz-Zygmund Inequality and in a special
context, but this will be enough for us. We followed
\cite[Chapter 10, Theorem 2, page 356]{Show-Teicher*97} by Y.S.~Show
and H.~Teicher. The relevant constant is the subject of
\cite{Ren-Liang*01} by Y.-F.~R and H.-Y.~Liang (their value is
slightly worse than ours) and \cite{Ferger*14} by D.~Ferger, where the
best constant is determined provided the $f_n$'s are ``symmetric''.

%%%%%%%%%%%%%%%%%% 
\begin{proof}
  We first notice that, since $\int_0^1 f_n(x)dx=0$, we may introduce
  a symmetrization through
  \begin{equation*}
    \sum_{n\le
      N}f_n(x_{2n-1})
    =-
    \int_{(x_{2n})\in X^{N}}
    \sum_{n\le
      2N}(-1)^nf_{\lceil{n/2\rceil}}(x_{n})d(x_{2n}).
  \end{equation*}
  Jensen's inequality gives us that
  \begin{multline*}
    \int_{(x_{2n-1})\in X^{N}}
    \biggl|
    \int_{(x_{2n})\in X^{N}}
    \sum_{n\le
      2N}(-1)^nf_{\lceil{n/2\rceil}}(x_{n})d(x_{2n})
    \biggr|^p
    d(x_{2n-1})
    \\\le
    \int_{(x_{2n-1})\in X^{N}}   
    \int_{(x_{2n})\in X^{N}}
    \biggl|\sum_{n\le
      2N}(-1)^nf_{\lceil{n/2\rceil}}(x_{n})\biggr|^p
    d(x_{2n})
    d(x_{2n-1})
   \end{multline*}
   from which we deduce that the $L^p$-norm of the symmetrization
   controls the one of the initial sum:
   \begin{equation*}
     \int_{(x_{2n-1})\in X^{N}}
     \biggl|\sum_{n\le
       N}f_n(x_{2n-1})
     \biggr|^p
    d(x_{2n-1})
     \le
    \int_{(x_{n})\in X^{2N}}
    \biggl|\sum_{n\le
      2N}(-1)^nf_{\lceil{n/2\rceil}}(x_{n})\biggr|^p
    d(x_{n}).
  \end{equation*}
  We next introduce Rademacher's system by noticing that, by
  successively exchanging $x_{2n-1}$ and $x_{2n}$, have 
  \begin{multline*}
    (1/2^N)\sum_{(\ve_n)\in\{\pm1\}^N}\int_{(x_{n})\in X^{2N}}
    \biggl|\sum_{n\le
      2N}\ve_{\lceil{n/2\rceil}}(-1)^nf_{\lceil{n/2\rceil}}(x_{n})\biggr|^pd(x_{n})
    \\=\int_{(x_{n})\in X^{2N}}
    \biggl|\sum_{n\le
      2N}(-1)^nf_{\lceil{n/2\rceil}}(x_{n})\biggr|^pd(x_{n}).
  \end{multline*}
  We may now remove the symmetrization since:
  \begin{align*}
    \int_{(x_{n})\in X^{2N}}
    &\biggl|\sum_{n\le 2N}\ve_{\lceil{n/2\rceil}}(-1)^nf_{\lceil{n/2\rceil}}(x_{n})\biggr|^p
    d(x_{n})
    \\&\le
      \int_{(x_{n})\in X^{2N}}
      2^{p-1}\biggl(
      \biggl|
      \sum_{n\le N}\ve_{n}f_n(x_{2n})\biggr|^p
      \biggl|
      +
      \sum_{n\le N}\ve_{n}f_n(x_{2n-1})\biggr|^p
      \biggr)
      d(x_{n})
    \\&\le
      2^p
      \int_{(x_{2n})\in X^{N}}
      \biggl|
      \sum_{n\le N}\ve_{n}f_n(x_{2n})\biggr|^p
      d(x_{2n}).
  \end{align*}
  The Khintchine Inequality from Theorem~\ref{KI} now gives us that
  \begin{equation*}
    (1/2^N)\sum_{(\ve_n)\in\{\pm1\}^N}\biggl|
      \sum_{n\le N}\ve_{n}f_n(x_{2n})\biggr|^p
    \le
    p^{p/2}\biggl(\sum_{n\le N}|f_n(x_{2n})|^2\biggr)^{p/2}.
  \end{equation*}
  The proof is then complete.
  Concerning the constant~4 in the Theorem, the paper
  \cite{Ren-Liang*01} by Y.-F.~Ren and H.-Y.~Liang gives the upper
  bound~9/2, which is worse than the above one.
\end{proof}
%%%%%%%%%%%%%%%%%%

%%%%%%%%%%%%%%%%%%%%%%%%%%%%%%%%%%%%%%%
\section{Proof of Theorem~\ref{CLS}}
%%%%%%%%%%%%%%%%%%%%%%%%%%%%%%%%%%%%%%%

%%%%%%%%%%
\begin{proof}
  We define $\Omega=\{1,\cdots, K\}$ which we equip with the
  probability measure defined by $\nu(\{k\})=|\lambda_k|/\|\lambda\|_1$.
  Given a positive integer $L$, we consider the family of functions
  $\varphi_\ell$, for $\ell\le L$ given by
  \begin{equation*}
    \begin{array}{rcl}
      \varphi_\ell:  \Omega^L\times X&\rightarrow&\mathbb{C}\\
      ((u_h)_{h\le L}, x)&\mapsto& \lambda_{u_\ell}^\circ g_{u_\ell}(x)
    \end{array}
  \end{equation*}
  so that
  \begin{equation*}
    \int_{\Omega^L}\varphi_\ell((u_h)_{h\le L}, x)
    d\nu = \sum_{k\in
      \Omega}\frac{|\lambda_k|}{\|\lambda\|_1}\lambda_{k}^\circ
    g_{k}(x) = \frac{f}{\|\lambda\|_1}=f_0
  \end{equation*}
  say. We aim at showing that $(1/L)\sum_{\ell\le
    L}\varphi_\ell((u_h)_{h\le L}, x)$ closely approximates~$f_0$ for
  most values of $(u_h)_{h\le L}$. Selecting one such value gives
  qualitatively our result. To do so, we write
  \begin{multline*}
    \int_{\Omega^L}\int_X
    \biggl|\frac{1}{L}\sum_{\ell\le
      L}\varphi_\ell((u_h)_{h\le L}, x)-f_0\biggr|^pd((u_h)_{h\le L})dx
    \\=
    \frac{1}{L^p}\int_{\Omega^L}\int_X
    \biggl|\sum_{\ell\le
      L}\bigl(\varphi_\ell((u_h)_{h\le L}, x)-f_0\bigr)\biggr|^pd((u_h)_{h\le L})dx.
  \end{multline*}
  We apply the Marcinkiewicz-Zygmung Inequality,
  i.e. Theorem~\ref{MZI}, to this latter expression, getting
  \begin{multline*}
    \int_{\Omega^L}\int_X
    \biggl|\frac{1}{L}\sum_{\ell\le
      L}\varphi_\ell((u_h)_{h\le L}, x)-f_0\biggr|^pd((u_h)_{h\le
      L})dx
    \\\le
    \frac{(4p)^{p/2}}{L^{p/2}}
    \int_{\Omega^L}\int_X
    \biggl|\frac{1}{L}\sum_{\ell\le
      L}\bigl|\varphi_\ell((u_h)_{h\le L}, x)-f_0\bigr|^2\biggr|^{p/2}d((u_h)_{h\le
      L})dx
    \\\le
    \frac{(4p)^{p/2}}{L^{p/2}}
    \int_X
    \int_{\Omega^L}\biggl|\frac{1}{L}\sum_{\ell\le
      L}\bigl|\varphi_1((u_h)_{h\le L}, x)-f_0\bigr|^2\biggr|^{p/2}d((u_h)_{h\le
      L})dx.
  \end{multline*}
  The end is straightforward:
  \begin{multline*}
    \int_{\Omega^L}\int_X
    \biggl|\frac{1}{L}\sum_{\ell\le
      L}\varphi_\ell((u_h)_{h\le L}, x)-f_0\biggr|^pd((u_h)_{h\le
      L})dx
    \\\le
    \frac{(4p)^{p/2}}{L^{p/2}}
    \int_X\int_{\Omega^L}\bigl|\varphi_1((u_h)_{h\le L}, x)-f_0\bigr|^{p}d((u_h)_{h\le
      L})dx.
  \end{multline*}
  Concerning the relevant $p$-norms, we make the following observations:
  \begin{multline*}
    \int_X\int_{\Omega^L}\bigl|\varphi_1((u_h)_{h\le L}, x)\bigr|^{p}d((u_h)_{h\le
      L})dx
    \\
    =\int_{\Omega^L}\biggl(
    \int_X\bigl|\varphi_1((u_h)_{h\le L}, x)\bigr|^{p}dx\biggr)d((u_h)_{h\le
      L})\le 1,
  \end{multline*}
  on the one side while on the other side, by the triangle inequality, we have
  \begin{equation*}
    \|f_0\|_p\le\sum_{k\le
      K}\frac{|\lambda_k|}{\|\lambda\|_1}\|g_k\|_p\le 1.
  \end{equation*}
  Therefore
  \begin{multline*}
    \biggl(\int_{\Omega^L}\int_X
    \biggl|\frac{1}{L}\sum_{\ell\le
      L}\varphi_\ell((u_h)_{h\le L}, x)-f_0\biggr|^pd((u_h)_{h\le
      L})dx\biggr)^{1/p}
    \\\le
    \frac{(4p)^{1/2}}{L^{1/2}}(1+1)=\sqrt{16p/L}.
  \end{multline*}
  We deduce from this inequality that the set of $(u_h)_{h\le L}$ for
  which
  \begin{equation*}
    \int_X
    \biggl|\frac{1}{L}\sum_{\ell\le
      L}\varphi_\ell((u_h)_{h\le L}, x)-f_0\biggr|^pdx > \ve^p
  \end{equation*}
  has measure at most $\sqrt{16p/(\ve^2 L)}$ which is strictly less
  than~1 by our assumption on~$L$. The theorem follows readily.
\end{proof}
%%%%%%%%%% 

%\bibliographystyle{authordate1}
\bibliographystyle{plain}
%\printbibliography
% $ biblatex auxiliary file $
% $ biblatex bbl format version 3.2 $
% Do not modify the above lines!
%
% This is an auxiliary file used by the 'biblatex' package.
% This file may safely be deleted. It will be recreated by
% biber as required.
%
% \bibliography{Local}

\end{document}